\documentclass[twoside]{article}
\def\nkyear{2018}
\def\nkvol{39}
\def\nkno{4}

\def\firstpage{1}
\def\lastpage{14}
\def\correction{4}
\def\shorttitle{Model Selection Consistency of Lasso for Empirical Data}
\def\shortauthor{{\it Y. Yang\ \ and\ \ H. Yang}}
\usepackage b \baselineskip 14.5pt
\usepackage{epstopdf}
\usepackage{multirow}
\usepackage[hang]{subfigure}
\usepackage[mathscr]{euscript}
\usepackage{srcltx}
\usepackage{bm}
\usepackage{amsmath}
\usepackage{amssymb}
\usepackage{amsthm}
\usepackage{graphicx}
\usepackage{color}
\usepackage{dsfont}

\newtheorem{rmk}{Remark}[section]

\newcommand\mR{\mathds {R}}

\newcommand\lam{\lambda}
\newcommand{\te}{\text{e}}
\newcommand{\e}{\epsilon}

\newcommand{\normm}[1]{\Big\vert\kern-0.25ex\Big\vert #1
		\Big\vert\kern-0.25ex\Big\vert}

\input{cam.doi}
\journalnumber{11401}
\DOImsnr{0001}
 \DOIyear{2018}

\title{\Large \bf \boldmath\ \\ Model Selection Consistency of Lasso for Empirical Data} 

\author{\large  Yuehan Yang$^1$\qquad  Hu Yang$^{2}$} 

\date{}

\begin{document}

\maketitle

\thispagestyle{first}
\renewcommand{\thefootnote}{\fnsymbol{footnote}}

\footnotetext{\hspace*{-5mm} \begin{tabular}{@{}r@{}p{13.4cm}@{}}
& Manuscript received November 8, 2016. Revised Feburary 24, 2017.  \\ 
$^1$ & School of Statistics and Mathematics, Central University of Finance and Economics, Beijing 100081, China.\\
&{E-mail:yyh@cufe.edu.cn} \\
$^{2}$ & College of Mathematics and Statistics, Chongqing University,
Chongqing 401331, China.\\
&{E-mail:yh@cqu.edu.cn}\\
 $^{\ast}$ & Project supported by
  National Natural Science Foundation of China (Grant No.11671059).
\end{tabular}}

\renewcommand{\thefootnote}{\arabic{footnote}}

\begin{abstract}
Large-scale empirical data, the sample size and the dimension are high,
often exhibit various characteristics. For example, the noise term follows unknown distributions or the model is very sparse that the number of critical variables is fixed while dimensionality grows with $n$. We consider the model selection problem of lasso for this kind of data. We investigate both theoretical guarantees and
simulations, and show that the lasso is robust for various kinds of data.

\vskip 4.5mm

\nd \begin{tabular}{@{}l@{ }p{10.1cm}} {\bf Keywords } &
Lasso, Model selection, Empirical data
\end{tabular}

\nd {\bf 2000 MR Subject Classification }
62J05, 62J07

\end{abstract}

\baselineskip 14pt

\setlength{\parindent}{1.5em}

\setcounter{section}{0}

\Section{Introduction}

Tibshirani \cite{Tibshirani1996} proposed the lasso (least absolute shrinkage and
selection operator) method for simultaneous model selection and estimation
of the regression parameters. It is very popular for high-dimensional
estimation due to its statistical accuracy for prediction and model
selection coupled with its computational feasibility. On the other hand,
under some sufficient condition the lasso solution is unique, and the
number of non-zero elements of lasso solution is always smaller than $n$ (see \cite{tib201201,wainwright2009sharp}). In recent years,
this kind of data has become more and more common in most fields. Similar
properties can also be seen in other penalized least squares since they
have a similar framework of solution.

Consider the problem of model selection in the sparse linear regression
model
\[
y_n = X_n \beta_n + \e_n,
\]
where the detail setting of the data can be found in the next section. Then the lasso estimator is defined as
\[
\widehat \beta_n(\lambda_n)\in  \mathop {\text{argmin}}\limits_{\beta\in \mR^p }\{\dfrac{1}{2}||y_n-X_n\beta||^2_2+\lambda_n ||\beta||\},
\]
where $\lam_n$ is the tuning parameter which controls the amount of
regularization. Set $
\widehat S_n\equiv\{j \in \{1,2,...,p_n\}: \widehat \beta_{j,n} \neq 0\}$ to select
predictors by lasso estimator $\widehat \beta_n$. Consequently, $\widehat S_n$ and
$\widehat \beta_n$ both depend on $\lam_n$, and the model selection criteria
results in the correct recovery of the set $S_n\equiv\{j \in \{1,2,...,p_n\}:  \beta_{j,n} \neq 0\}$:
\[
P(\widehat S_n = S_n) \rightarrow 1 \ \text{as} \ n \rightarrow \infty.
\]

On the model selection front of the lasso estimator, Zhao and Yu \cite{YuBin06(lasso)}
established the irrepresentable condition on the generating covariance
matrices for the lasso's model selection consistency. This condition was
also discovered in \cite{meinshausen2006high,yuan2006model,zou2006adaptive}. Using
the language of \cite{YuBin06(lasso)}, the irrepresentable condition is defined
as $
|C_{21}C^{-1}_{11}sign(\beta_{(1)})|\leqslant \mathbf{1} -\eta$, where
sign($\cdot$) maps positive entry to $1$, negative entry to $-1$ and zero to
zero. The definitions of $C_{21}$ and $C_{11}$ can be seen in Section 2.
When signs of the true coefficients are unknown, they need $l_1$ norms of
the regression coefficients to be smaller than $1$. Beyond lasso, regularization methods also have been widely used for
high-dimensional model selection, e.g.,
\cite{bickel2009simultaneous,candes2007dantzig,fan2001variable,Fan2004,lv2009unified,meinshausen2010stability,yyh(elastic),yyh(lasso),yyh(adaptive),zhang2010nearly,zou2005elastic,zou2009adaptive}.
There has been a considerable amount of recent work dedicated to the lasso
problem and regularization methods problem.

Yet, the study of model selection problem for empirical data is still needed. Stock data for instance, the Gaussian assumption of the noise term is always unsatisfied for these data. And the critical variables are extremely few contrast to the collected dimensionality. In this paper, we consider this kind of data: The sample size and the dimension are high, but the information of critical variable data is missing (The signs of the true $ \beta_n $ and the distribution of the noise terms are unknown.) and the model is extremely sparse that the number of nonzero coefficients is fixed. This kind of data are common in the empirical analysis hence we called it empirical data.

We consider the model selection consistency of lasso and investigate regular conditions to fit the data setting. Under conditions, the probability for lasso to select the true model is covered by the probability of
\[
\{ ||W_n||_{\infty}  \leqslant G_n\},
\]
where $W_n = X'_n \e_n / \sqrt{n}$ and $ G_n $ is a function of $ \lam_n $,
$ n $, $ q $. Above inequality is simple and also easy to calculate its
probability. Based on the train of thought of the proof, we analyze the
model selection consistency of lasso under easier conditions than the irrepresentable condition for empirical data. In the simulation part, we discuss the effectiveness of lasso. Four samples are given, in which the irrepresentable condition fails for all the settings, but lasso still can select variables correctly in two of them when our conditions hold.

We discuss the different assumptions of noise terms
$\epsilon_{i,n}$ for model selection consistency. Gaussian errors or the
subgaussian errors\footnote{e.g. $P(|\e_{i,n}| \geqslant t) \leqslant C
\te^{-ct^2}, \ \forall t \geqslant 0$.} would be standard, but possess a strong
tail. One basic assumption in this paper is that, errors are assumed to be identically and independently distributed with zero mean and finite variance.

The rest of the paper is organized as follow. In Section 2, we investigate the data setting, notations, and conditions. We introduce a lower bound to cover the case in which the lasso chooses wrong models when suitable conditions hold. Then, to demonstrate the advantages of this bound, we show the different settings and different assumptions of noise terms in Section 3. We show that the lasso has model selection consistency for empirical data with mild conditions. Section 4 presents the results of the simulation studies. Finally, in Section 5, we present the proof of the main theorem.

\Section{Data Setting, Notations and Conditions}

Consider the problem of model selection for specific data
\[
y_n = X_n \beta_n + \e_n,
\]
where $\epsilon_n = (\epsilon_{1,n},\epsilon_{2,n},...,\epsilon_{n,n})'$ is a vector of i.i.d. random variables with mean $0$ and variance $\sigma^2$, $X_n $ is an $n\times p_n$ design matrix of predictor variables, $\beta_n \in \mathds {R}^{p_n}$ is a vector of true regression coefficients and is commonly
imposed to be sparse with only a small proportion of nonzeros. Without loss
of generality, write $\beta_n=(\beta_{1,n},...,\beta_{q,n},
\beta_{q+1,n},...,\beta_{p_n,n})'$ where $\beta_{j,n} \neq 0$ for $j=1,...,q$
and $\beta_{j,n}=0$ for $j =q+1,...,p_n$. Then write
$\beta^{(1)}_n=(\beta_{1,n},...,\beta_{q,n})'$ and
$\beta^{(2)}_n=(\beta_{q+1,n},...,\beta_{p_n,n})$, that is, only the first $q$
entries are nonvanishing. Besides, for any vector $ \alpha=
(\alpha_1,...,\alpha_m)' $, we denote $ ||\alpha||  = \sum_{i=1}^{m}
|\alpha_i|$, $ ||\alpha||_2^2 = \sum_{i=1}^{m} \alpha_i^2 $ and $
||\alpha||_{\infty} = \max\limits_{i=1,...,m} |\alpha_i| $.

For deriving the theoretical results, we write $X_{n}(1)$ and $X_n(2)$ as
the first $q$ and the last $p_n-q$ columns of $X_n$,S respectively. Let $C_n =
\dfrac{1}{n} X'_nX_n$. Partition $C_n$ as
\[
C_n = \left( {\begin{array}{*{20}c}
	{C_{11,n}} \quad {C_{12,n}}\\
	{C_{21,n}} \quad {C_{22,n}} \\
	\end{array}} \right),
\]
where $C_{11,n}$ is $q \times q$ matrix and assumed to be invertible. Set
$W_n=X'_n \e_n/ \sqrt{n}$. Similarly, $W^{(1)}_n$ and $W^{(2)}_n$ indicate
the first $q$ and the last $p_n-q$ elements of $W_n$. Suppose that $\Lambda_{min}(C_{11,n}) >0$ denotes the smallest eigenvalue of
$C_{11,n}$ and consider that $ q $ does not grow with $ n $. We introduce the following conditions.

\begin{enumerate}
	\item[(C1)] For $j=q+1,...,p_n$, let $e_j$  be the unit vector in the direction of $ j $-th coordinate. There exists a
	positive constant $ 0<\eta<1 $ such that
	\[
	||e'_jC_{21}||_2 \leqslant (1-\eta).
	\]
\end{enumerate}

\begin{enumerate}
	\item[(C2)]  There exists $\delta \in (0,1)$, such that for all $n >
	\delta^{-1}$ and $x \in \mR^q$, $y \in \mR ^{p_n-q}$,
	\[
	(x'C_{12,n}y)^2 \leqslant \delta^2 (x'C_{11,n}x)\cdot(y'C_{22,n}y).
	\]
\end{enumerate}

(C1) and (C2) play a central role in our theoretical analysis. Both conditions are easy to satisfy. (C1) for instance, it requires an upper bound on $l_2$-norm, which is much weaker than requires the upper bound on $l_1$-norm, i.e. irrepresentable condition and variants of this condition
\cite{fan2014adaptive,huang2008asymptotic,meinshausen2006high,YuBin06(lasso),zou2009adaptive}. Another advantage of (C1) is that we do not need the signs of the true coefficients.
(C2) requires that the multiple correlation between relevant variables and the irrelevant variables is strictly less than one. It is weaker than assuming orthogonality of the two sets of variables. This condition also has regular appeared many times in the literature, for example, \cite{meinshausen2009lasso}.

Then we have the following theorem, which describes the relationship between the probability
of lasso choosing the true model and the probability of $\{
||W_n||_{\infty}  \leqslant G_n \}$. Videlicet, it is a lower bound on the
probability of lasso picking the true model.

\begin{thm}
	\label{thm1}
	Assume that (C1)-(C2) hold and $0< \Lambda_{min} \triangleq  \Lambda_{min}(C_{11,n}) $. Set $ \rho
	\in (0,1) $. We have
	\[
	P(\widehat S_n = S_n) \geqslant P(||W_n||_{\infty} \leqslant G_n),
	\]
	where $ G_n = \min \left\{
	\sqrt{n} \Lambda_{min}\left(
	\min\limits_{j = 1,...,q}|\beta_{j,n}|
	-\dfrac{\lam_n\sqrt{q}}{\Lambda_{min}\cdot n} \right),
	\ \dfrac{ \lam_n \rho}{\sqrt{n}} \right\} $.
\end{thm}

\begin{rmk}
	Theorem~\ref{thm1} is a key technical tool in the theoretical results. It puts a
	lower bound on the probability of lasso selecting the true model, and
	this bound is intuitive to calculate. Besides that, considering about $
	G_n $, it is easy to find out that there exists a lower bound of non-zero
	coefficients $ \min\limits_{j = 1,...,q}| \beta_{j,n}|  >
	\lam_n\sqrt{q}/(\Lambda_{min}\cdot n)$. This bound can be controlled by
	the regularization parameter $\lam_n$. It is also a regular assumption
	in the literature that the non-zero coefficients cannot be too small.
\end{rmk}
\begin{rmk}
	According to the proof of Theorem~\ref{thm1}, we can find that it is
	also directly to obtain the sign consistency of the lasso (see the
	latter part of the proof). Besides, Theorem~\ref{thm1} can be applied in
	a wide range of dimensional setting. We will discuss the behavior of the
	lasso on model selection consistency under different settings in the
	next section.
\end{rmk}

\Section{Model selection consistency}

Now we consider the decay rate of the probability of $\{||W_n||_{\infty} >
G_n \}$. Different dimensions and different assumption of noise terms are
discussed in this section.

First, we consider general dimensional setting, i.e. $p_n= O(n^{c_1})$ where $0 < c_1 <1$. Under this setting, we can obtain
the model selection consistency of lasso by no constraint for the noise
terms. Then, we consider ultra-high dimensional
setting, i.e. $p_n = O(\te^{n^{c_2}})$ where $ 0 < c_2 <1 $. Under this setting, we need an assumption of $\e_n$ to make the
model selection of lasso. Gaussian assumption would be a simple and common
one, but result in a strong tail. We prefer the more standard assumption
that only i.i.d random variables of the noise terms.

Before discussing the detail rate of the probability of lower bound, we give
the following regular condition
\begin{enumerate}
	\item[(C3)]  $n^{-1} X'_{j,n}X_{j,n} \leqslant 1$ for $j=1,...,p_n$.
\end{enumerate}

It is a typical assumption in sparse linear regression literature. It can be
achieved by normalizing the covariates. See
\cite{huang2008asymptotic,YuBin06(lasso)}.

\Subsection{General dimensional setting $p_n = O(n^{c_1})$}

In this part, we consider the general dimensional setting where $p_n$ is allowed to grow with $n$ and show the model selection consistency of lasso as follows
\begin{thm}
	\label{thm2}
	Assume that $\e_i$ are i.i.d random variables with mean $0$ and variance
	$\sigma^2$. Suppose (C1)-(C3) hold. For $p_n = O(n^{c_1})$ where $ 0<c_1<1 $, if $ \frac{\lam_n}{\sqrt{n}} \propto n^{c_3/2}$ where $ c_1 < c_3 < 1 $ and $ \min\limits_{j = 1,...,q}|
	\beta_{j,n}| >n^{\frac{c_3-1}{2}} $, then we have
	\[
	P(\widehat S_n = S_n) \geqslant  1- n^{c_1-c_3} \rightarrow 1 \ \text{as} \
	n\rightarrow \infty.
	\]
\end{thm}

\noindent\textbf{Proof,}

Following the result in Theorem \ref{thm1}, we have
\[
P(\widehat S_n = S_n) \geqslant P(||W_n||_{\infty} \leqslant G_n),
\]
where
\[
G_n = \min \left\{ \sqrt{n} \Lambda_{min}\left(\min\limits_{j = 1,...,q}|
\beta_{j,n}|  -\dfrac{\lam_n\sqrt{q}}{\Lambda_{min}\cdot n}  \right), \
\dfrac{ \lam_n \rho}{\sqrt{n}}  \right\}.
\]

Applying the setting of Theorem~\ref{thm2}, hence for $ n \rightarrow \infty $,
\[
\Lambda^{-1}_{min} \cdot \lam_n\dfrac{\sqrt{q}}{n}=O(n^{\frac{c_3-1}{2}})
\rightarrow 0.
\]

Then there exists a positive constant $ K_n $ that
\[
G_n = \rho \lam_n / \sqrt{n}= K_n n^{c_3 / 2}.
\]

If (C3) holds, by Markov's inequality, we easily get
\begin{align}
P(||W_n||_{\infty} > G_n ) &\leqslant \sum^{p_n}_{j=1} P(|W_{j,n}| > G_n)
\nonumber\\
& = \sum^{p_n}_{j=1}  P\Big(\Big|\dfrac{X'_{j,n} \e}{\sqrt{n}}\Big|> K_n n^{c_3/2} \Big) \nonumber\\
& \leqslant K^{-2}_n n^{-c_3} \cdot n^{c_1}  \rightarrow 0 \ \text{as} \ n
\rightarrow \infty. \nonumber
\end{align}
The proof is completed.
\\

The proof of Theorem~\ref{thm2} states that in this setting, lasso is robust and selects the true model with regular restraints. Similarly, if we consider the classical setting
where $p$, $q$ and $\beta$ are fixed when $n\rightarrow \infty$, we have the following result
\begin{coro}\label{coro}
For fixed $p$, $q$ and $\beta$, under regularity assumptions (C1)-(C3), assume $\e_i$ are i.i.d random variables with mean $0$ and variance
$\sigma^2$. If $ \lam_n $ satisfies that
$\frac{\lam_n}{\sqrt{n}} \rightarrow \infty$ and $ \frac{\lam_n}{n}
\rightarrow 0 $ when $ n \rightarrow \infty $, then
\[
P(\widehat S = S)  \rightarrow 1 \ \text{as} \ n\rightarrow \infty.
\]
\end{coro}
Similar with the argument of Theorem~3.1, Corollary \ref{coro} can be proved directly by Markov's inequality, hence the proof is omitted here.

Besides, if we assume that the noise term follows the Gaussian assumption, under the same setting of Theorem~3.1, we have
\begin{align}
\label{eq:S1}
P(\widehat S_n \neq S_n) & \leqslant P(||W_n||_{\infty} > G_n)  \nonumber\\
&\leqslant \sum^{p_n}_{j=1} P(|W_{j,n}| > K_n n^{c_3/2}) \nonumber\\
& <  n^{ c_1-c_3 /2} \te^{-\frac{1}{2}n^{c_3}}
\rightarrow 0 \ \text{as} \ n \rightarrow \infty,
\end{align}

\Subsection{Ultra-high dimensional setting $p_n = O(\te^{n^{c_2}})$}

In this part, we consider the ultra-high dimensional setting as $p_n =
O(\te^{n^{c_2}})$ where $0<c_2<1$ and discuss the different situation
under different assumptions of noise terms (Gaussian assumption and
non-Gaussian assumption). Theorem~\ref{thm3} shows the result under
non-Gaussian assumption by applying Bernstein's inequality. Also, we show
the model selection consistency of the lasso under the Gaussian assumption in Corollary~\ref{coro2}.

We shall make use of the following condition
\begin{enumerate}
	\item[(C4)] Assume that $\e_{1,n},...,\e_{n,n}$ are independent random variables with mean $0$ and the following inequality satisfies for $j=1,...,p_n$
	\[
	\dfrac{1}{\sqrt{n}}E|W_{j,n}|^m \leqslant \dfrac{m!}{2}L^{m-2}, m = 2,3,...
	\]
	where $W_{j,n} = \dfrac{1}{\sqrt{n}} X'_{j,n} \e_n$ and $L \in (0, \infty)$.
\end{enumerate}

(C4) is the precondition for the non-Gaussian assumption (The model selection consistency of the lasso under the Gaussian assumption does not need this condition.). It is applied here for the Bernstein's inequality. According to (C4), we have
\[
E \exp\left[X'_{j,n} \cdot \e_n/ L_0 \right]\leqslant \exp\left[
\dfrac{n}{2(L^2_0 - L \cdot L_0)}
\right],
\]
where $ L_0 >L $. This bound leads to Bernstein's inequality as given in \cite{bennet1962probability}. Then we have the following result.

\begin{thm}
	\label{thm3}
	Assume that $\e_i$ are i.i.d random variables with mean $0$ and variance $\sigma^2$, suppose (C1)-(C2) and (C4) hold. Set $p_n = O(\te^{n^{c_2}})$ where $ 0 < c_2 <1 $. If $\dfrac{\lambda_n}{\sqrt{n}} \propto n^{c_3/2}$ where $c_2< c_3<1$ and $\min\limits_{j=1,...,q}|\beta_{j,n}| > n^{\frac{c_3-1}{2}}$. We have
	\[
	P(\widehat S_n = S_n) \geqslant 1- \te^{-n^{c_2}}\rightarrow 1 \ \text{as} \ n\rightarrow \infty.
	\]
\end{thm}

\noindent\textbf{Proof:}

By Bernstein's inequality, let $ t > 0 $ be arbitrary. We have
\[
P(W_{j,n} \geqslant n^{c_3/2}(Lt + \sqrt{2t})) \leqslant \te^{-tn^{c_3}} \leqslant \te^{-n^{c_2}} .
\]

Applying the result of Lemma~14.13 from \cite{buhlmann2011high}, when (C3) holds, we have
\[
 P(\max\limits_{1 \leqslant j \leqslant p_n}|W_{j,n}| \geqslant n^{c_3/2}(Lt + \sqrt{2t} + \alpha(L,n,p_n))) \leqslant  \te^{-n^{c_2}}.\]

Following the setting of $p$,
\[
\alpha(L,n,p_n)  = \sqrt{(2\log2p_n)/n}+ (L\log(2p_n))/n \rightarrow 0.
\]

Let $J \in (0, \infty)$ to make following inequalities hold for all $t >0$,
\[
\alpha(L,n,p_n) < J, \ Lt + \sqrt{2t} \leqslant Jt.
\]

Then we have
\[
P(||W_n||_{\infty} > J (1+t) n^{c_3/2}) \leqslant  \te^{-n^{c_2}},
\]
which complets the proof.
\\

Similarly as in general high-dimensional setting, we have the following result under Gaussian assumption. Since the proof of Corollary~\ref{coro2} is direct, we just state the result here without proof.

\begin{coro}
	\label{coro2}
	Assume that $ \e_i  $ are i.i.d Gaussian random variables. Let $ p_n = O(\te^{n^{c_2}}) $ where $ 0 < c_2 < 1 $.
	Suppose (C1) - (C3) hold. If $\frac{\lam_n}{\sqrt{n}}  =O(n^{c_3/2})$ where $c_2 < c_3 < 1$ and $ \min\limits_{j = 1,...,q}| \beta_{j,n}| >n^{\frac{c_3-1}{2}} $. The Gaussian assumption of noise terms are considered in the following
	\begin{align}
	P(\widehat S_n \neq S_n) & \leqslant P(||W_n||_{\infty} >G_n )  \nonumber\\
	&\leqslant \sum^{p_n}_{j=1} P(|W_{j,n}| > G_n ) \nonumber\\
	& = O( n^{-c_3/2} \te^{n^{c_2}-\frac{1}{2}n^{c_3}}) = o(\te^{-n^{c_2}}) \rightarrow 0 \ \text{as} \ n \rightarrow \infty. \nonumber
	\end{align}
\end{coro}

\Section{Simulation Part}

In this section, we evaluate the finite sample property of lasso estimator with synthetic data. We start with the behavior of lasso under different settings, then considering the relationship between $n$, $p$, $q$ and then consider the different noise terms.

\Subsection{Model selection}

This first part illustrates two simple cases (low dimension vs high dimension) to show the efficiency of lasso. Following cases describe two different settings to lead the lasso's model selection consistency and inconsistency when (C1) and (C2) hold and fail. As a contrast, we introduce the irrepresentable condition in this part, and it fails in all the settings.

\textbf{Example. 1}

In the low dimensional case, assume that there are $n=100$ observations and the values of parameters are chosen as $p = 3$, $q = 2$, that is,
\[
\beta = \{2,3,0\}.
\]

We generate the response $y$ by
\[
y=X_1\beta_1+X_2 \beta_2+ X_3\beta_3+\epsilon,
\]
where $X_1$, $X_2$ and $\e$ are i.i.d random variables from Gaussian distribution with mean $0$ and variance $1$. The third predictor $X_3$ is generated to be correlated with other parameters as the following two cases:
\[
X_3 = \dfrac{2}{3} X_1 + \dfrac{2}{3} X_2 + \dfrac{1}{3}e
\]
and
\[
X_3 = \dfrac{1}{2} X_1 + \dfrac{1}{2} X_2 + \dfrac{1}{\sqrt{2}}e,
\]
where $e$ is i.i.d random variable with the same setting as $\e$.

\begin{figure}[!ht]
	\centering
	\subfigure[Lasso fails]{\label{fig:ec1}
		\includegraphics[width=0.48\textwidth]{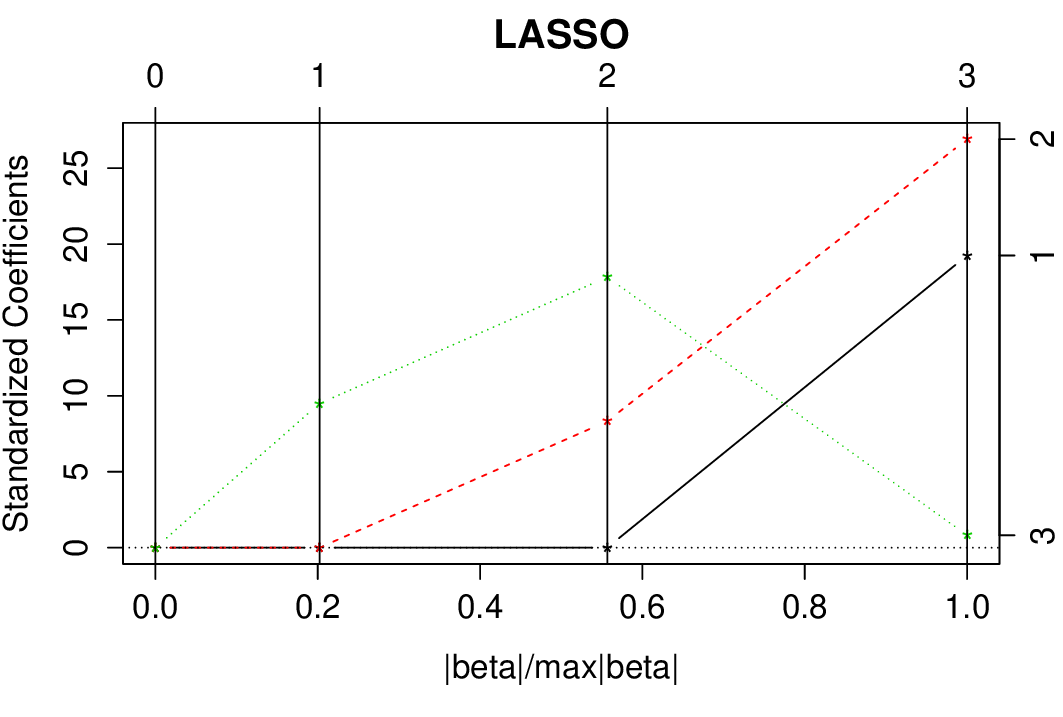}}
	\subfigure[Lasso holds]{\label{fig:ec2}
		\includegraphics[width=0.48\textwidth]{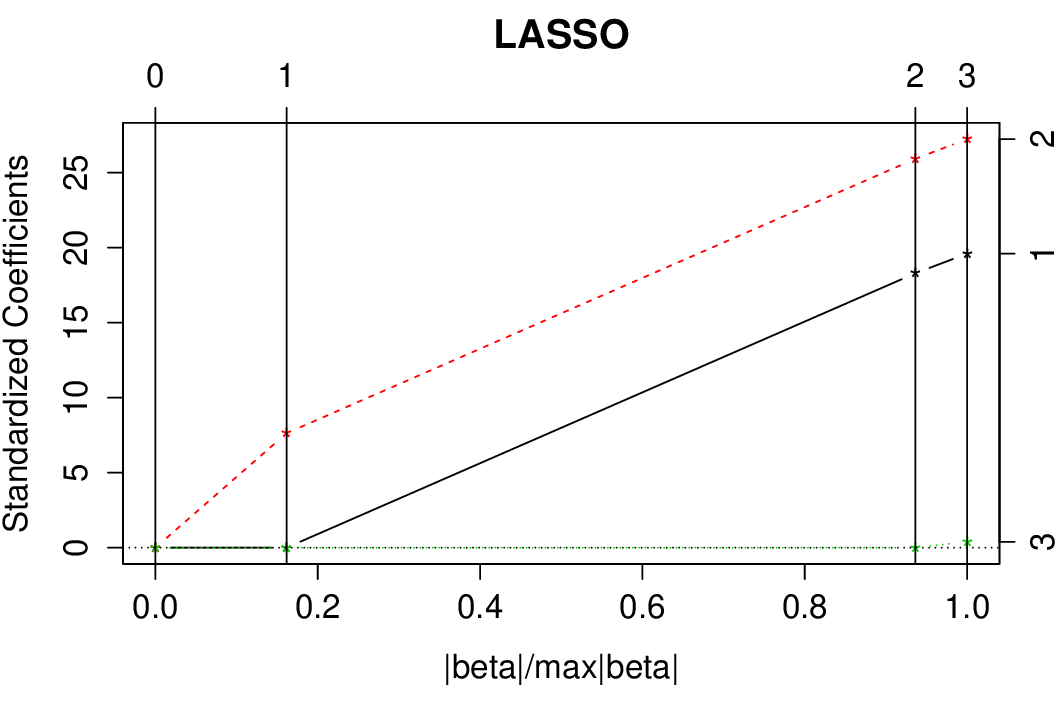}}
	\subfigure[Lasso fails]{\label{fig:ec3}
		\includegraphics[width=0.48\textwidth]{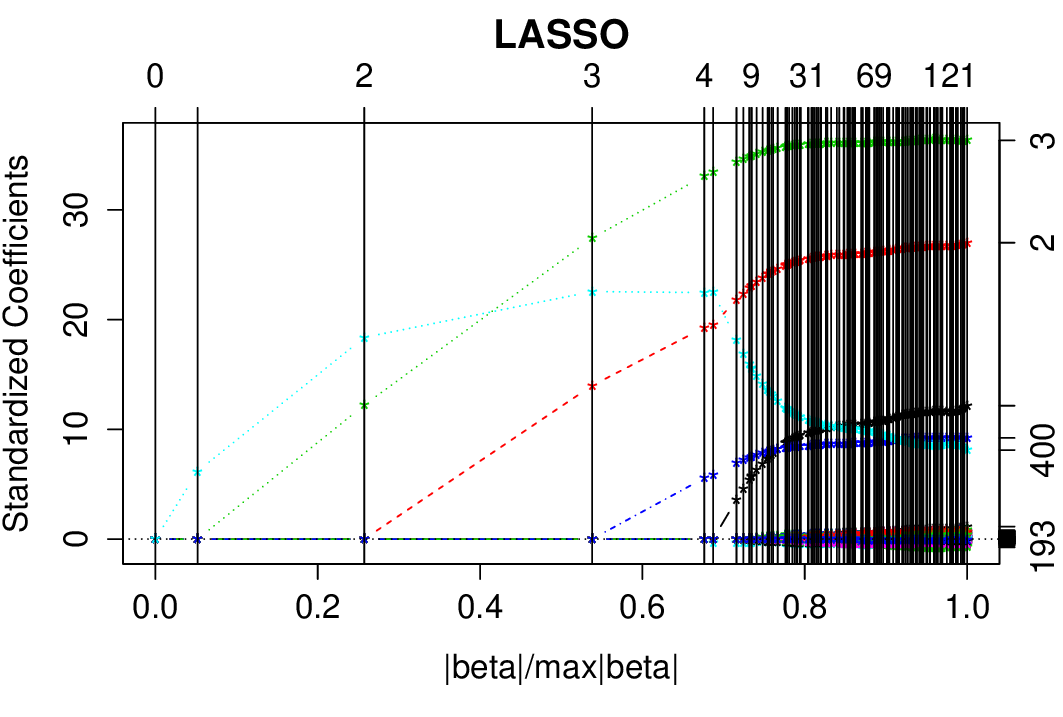}}
	\subfigure[Lasso holds]{\label{fig:ec4}
		\includegraphics[width=0.48\textwidth]{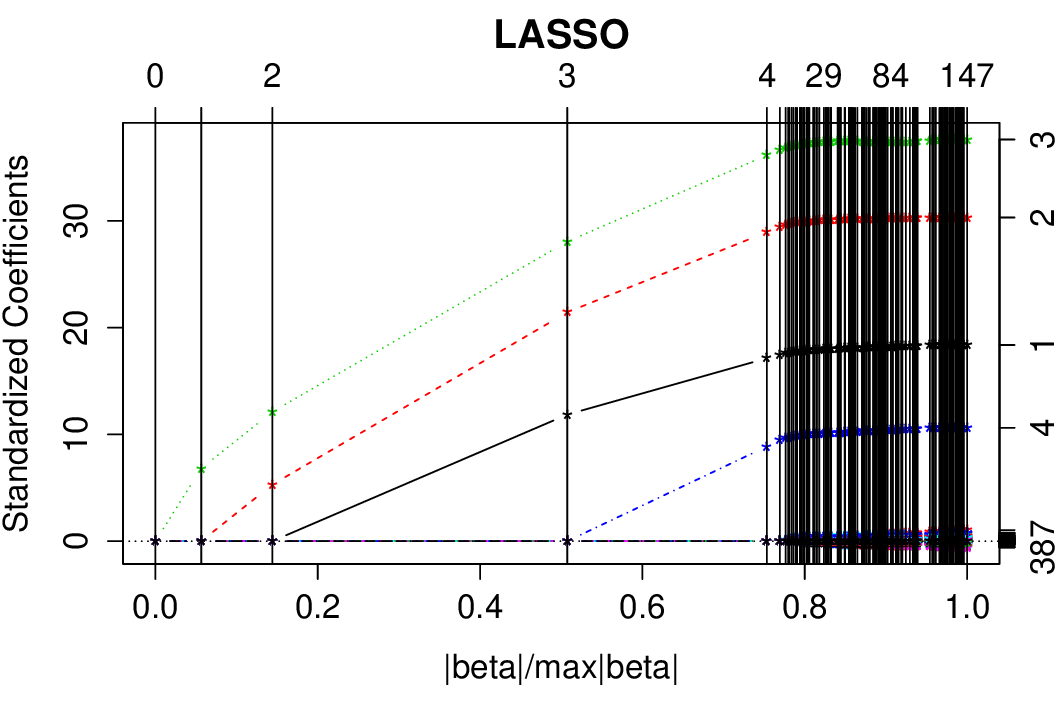}}
	\caption{An example to illustrate the efficiency of lasso's (in)consistency in model selection. Above two graphs are constructed in a low dimensional setting. Below graphs are constructed in a high dimensional setting. The left graphs in both columns are set where (C1) and (C2) fail, and the right graphs are set where (C1) and (C2) hold.}
	\label{fig:1}
\end{figure}

We can find that the lasso fails for the first case when (C1) and (C2) fail, and selects the right model  for the second case when (C1) and (C2) hold. The different solutions are illustrated by Figure \ref{fig:1}. Since the irrepresentable condition fails in both cases, it shows that the lasso suits more kinds of data even if the irrepresentable condition is relaxed.

\textbf{Example. 2}

We construct a high dimensional case which $p=400$, $q=4$ and $n=100$. The true parameters are set as
\[
\beta = \{2,3,1,4,0,0,...,0\}
\]
and the response $y$ is generated by
\[
y=X\beta+\epsilon,
\]
where $X = (X_1,...,X_p)$ is $100 \times 400 $ matrix, and the elements of $X$ are i.i.d random variables from Gaussian distribution with mean $0$ and variance $1$ except $X_{400}$. The last predictor $X_{400}$ is generated respectively as follows in the two settings
\[
X_{400} = \dfrac{7}{8} X_1 + \dfrac{3}{8} X_2+\dfrac{1}{8} X_3
+\dfrac{1}{8} X_4+\dfrac{1}{8} X_5+\dfrac{1}{8} X_6+\dfrac{1}{8} X_7+ \dfrac{1}{8}e
\]
and
\[
X_{400} = \dfrac{1}{4} X_1 + \dfrac{1}{4} X_2+\dfrac{1}{4} X_3
+\dfrac{1}{4} X_4+\dfrac{1}{4} X_5+\dfrac{1}{4} X_6+\dfrac{1}{4} X_7+ \dfrac{3}{4}e,
\]
where $e$ follows the same setting as Example 1. Hence $X_{400}$ is also constructed from Gaussian distribution with mean $0$ and variance $1$. We find that our conditions also fail for the first high dimensional case but hold for the second. Besides that, irrepresentable condition fails for both two situations.

We get different lasso solutions for above four cases in Figure \ref{fig:1} (the lasso path is got by lars algorithm \cite{Efron}).
As shown in Figure \ref{fig:1}, both graphs on the left satisfy neither irrepresentable condition nor (C1)-(C2), and lasso cannot select variables correctly (Both graphs select other irrelevant variables, e.g., $X_4$ in the first graph and $X_{400}$ in the second.). In contrast, both graphs on the right select the right model in the settings that (C1) and (C2) hold and irrepresentable condition fails.

Besides that, the above examples are all constructed based on the synthetic data, in which the unknown parameter is actually known. In the empirical analysis, the true model cannot be known in advance. We should recognize a situation in which lasso can be used without precondition.

\Subsection{Relationship between $p$, $q$ and $n$}

In this part, we give a direct view to show the relationship between $n$, $p$ and $q$, or to say how the sparsity and the sample size affect the model selection of lasso.

The nonzero elements $\beta^{(1)}$ are set as
\[
\{9,6,8,12,19,8,19,9,6,8,12,19,8,19\}.
\]

If the number of nonzero elements is less than $14$, we select the number in sequence. The rest of the other elements in this gather are shrunk to zero. The number of observations and the parameters are chosen as Table \ref{table:1}. The predictors are made from Gaussian random generation. Among this table, lasso selects the right variables in the first six items in the list and selects the wrong variables in the remaining items in the list.

The high dimensional settings are considered. The results indicate that $q$ is always required to be small enough for the efficiency of the lasso. When the number of critical factors increases, the sample size needs to be increased too to make sure the lasso chooses the right model. In contrast, the number of zero elements has less influence on the lasso's (in)consistency in model selection.

\begin{table}[!ht]
	\centering
	\caption{Example settings}
	\label{table:1}
	\begin{tabular}{cccc|ccccc}
		\hline
		Example	& $n$   & $p$ & $q$&	Example	& $n$   & $p$ & $q$&\\
		\hline
		1& 100 & 400& 4 &  7& 100 & 400& 5 &  \\
		\hline
		2&  100  &  500& 5 & 8& 100 & 500& 6 &\\
		\hline
		3	&  200 &  500 & 7  & 9& 100 & 500& 7&\\
		\hline
		4	&    200   &  1000&  7 & 10 & 100 & 1000& 7&\\
		\hline
		5	&  500  &  500 & 14 & 11 & 100 & 2000& 7&\\
		\hline
		6	& 500 & 2000  &  14 & 12 & 300 & 2000& 14&\\
		\hline
	\end{tabular}
\end{table}

\Subsection{Different noise terms}

In this part, we consider a high dimensional example with different noise terms. Data from the high-dimensional linear regression model is set as
\[
y_i= X'_i \beta + \epsilon_i, \ \ i=1,...,n,
\]
where the data has $n=100$ observations and the value of parameter is chosen as $p = 1000$. The true regression coefficient vector is fixed as
\[\beta= \{9,6,8,0,...,0\}.\]

For the distribution of the noise $\epsilon$, we consider four distributions: Gaussian assumption with mean $0$ and variance $1$; exponential distribution with rate $1$; uniform distribution with minimum $0$ and maximum $1$; student's t with degrees of freedom $100$.

\begin{figure}[!ht]
	\centering
	\subfigure[Gaussian assumption]{\label{fig:11}
		\includegraphics[width=0.48\textwidth]{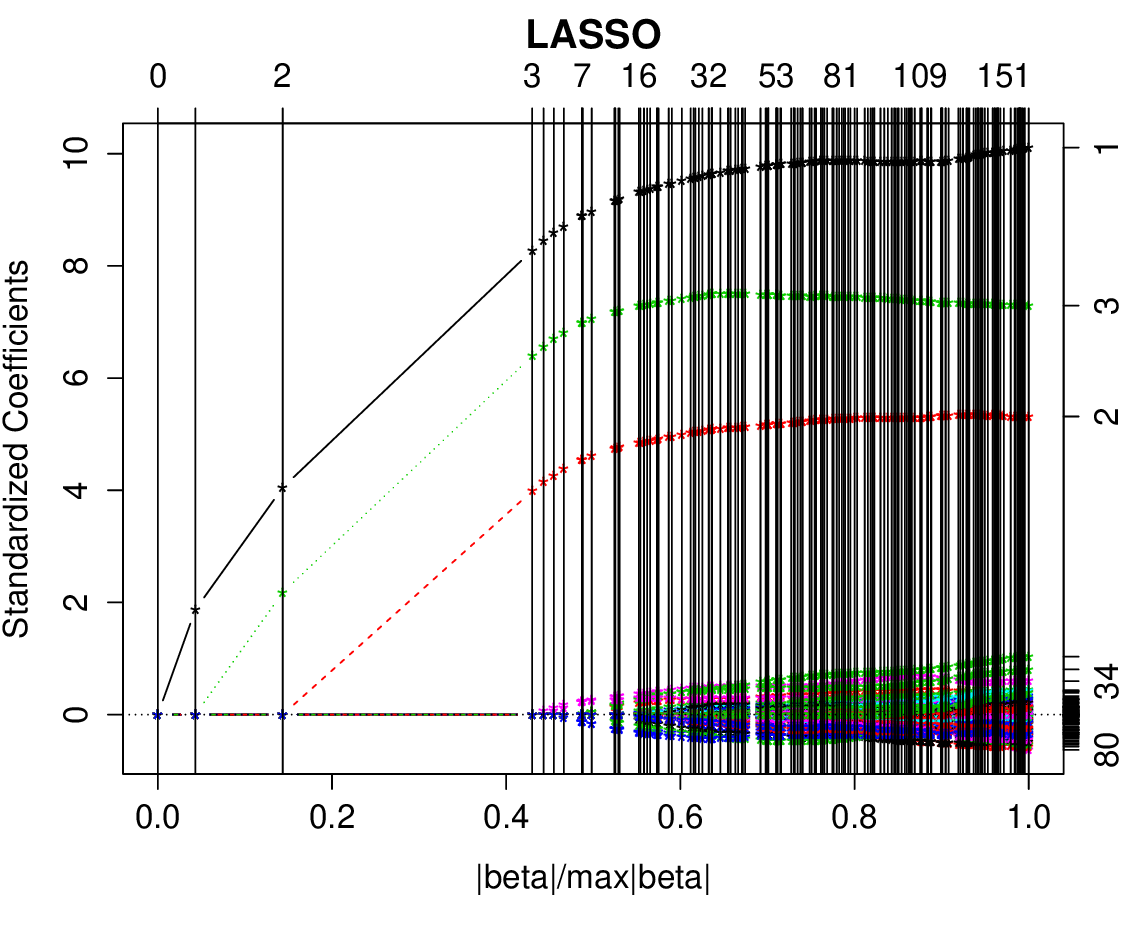}}
	\subfigure[Exponential distribution]{\label{fig:12}
		\includegraphics[width=0.48\textwidth]{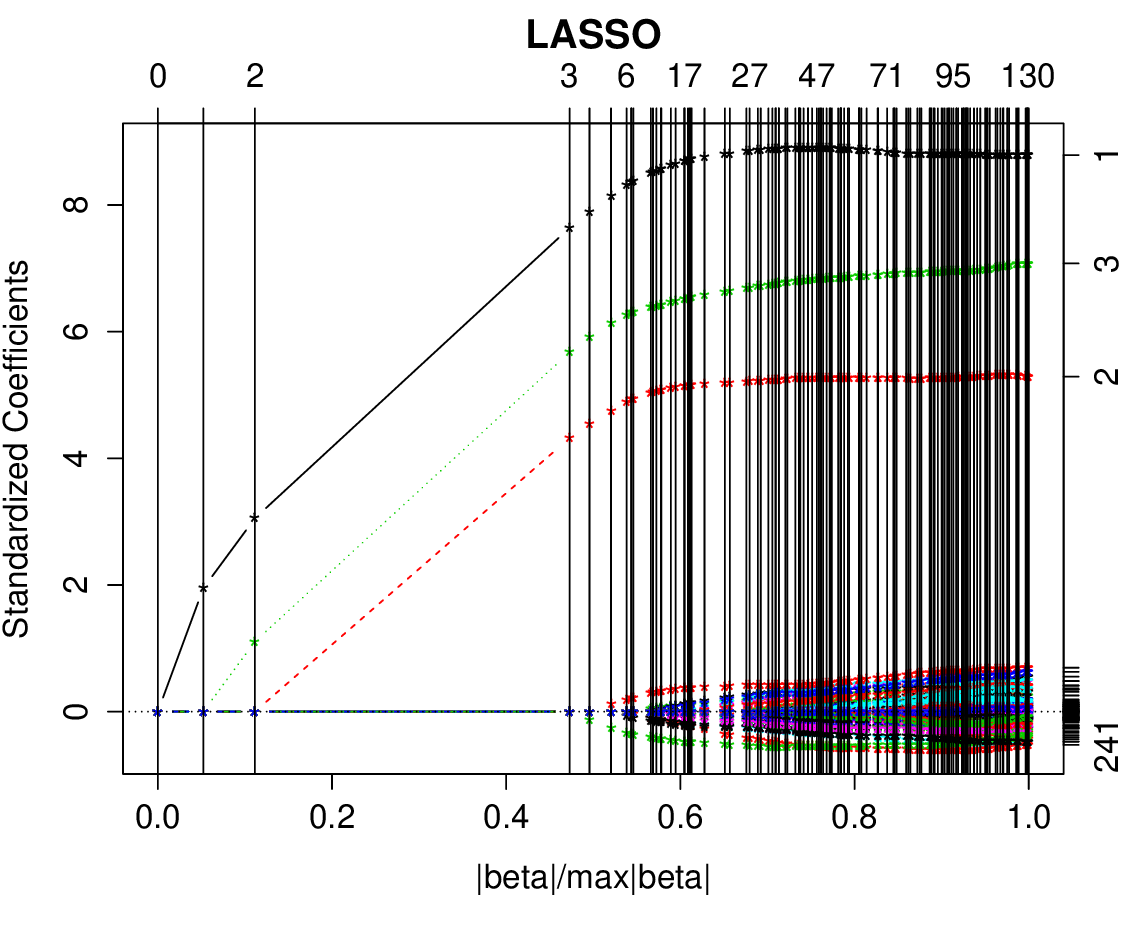}}
	\subfigure[Uniform distribution]{\label{fig:13}
		\includegraphics[width=0.48\textwidth]{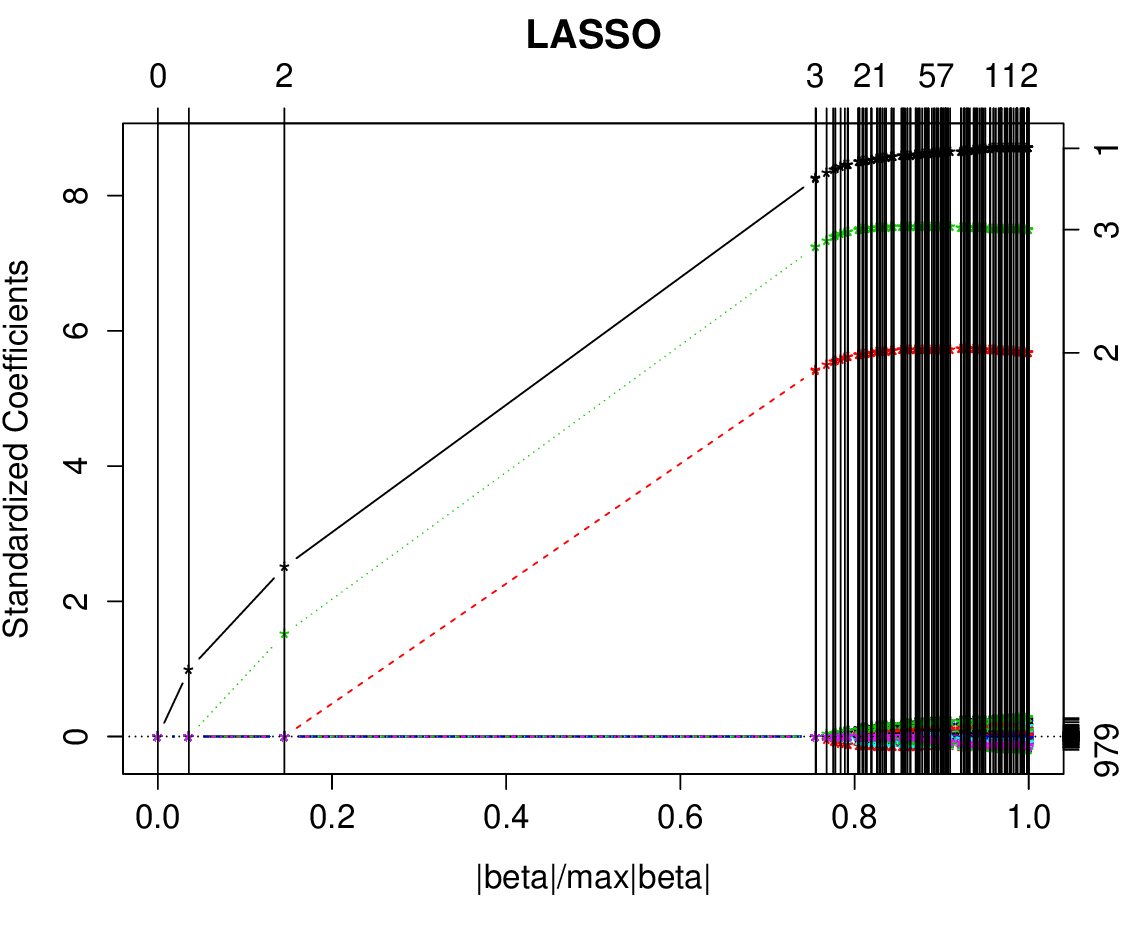}}
	\subfigure[Student's t]{\label{fig:14}
		\includegraphics[width=0.48\textwidth]{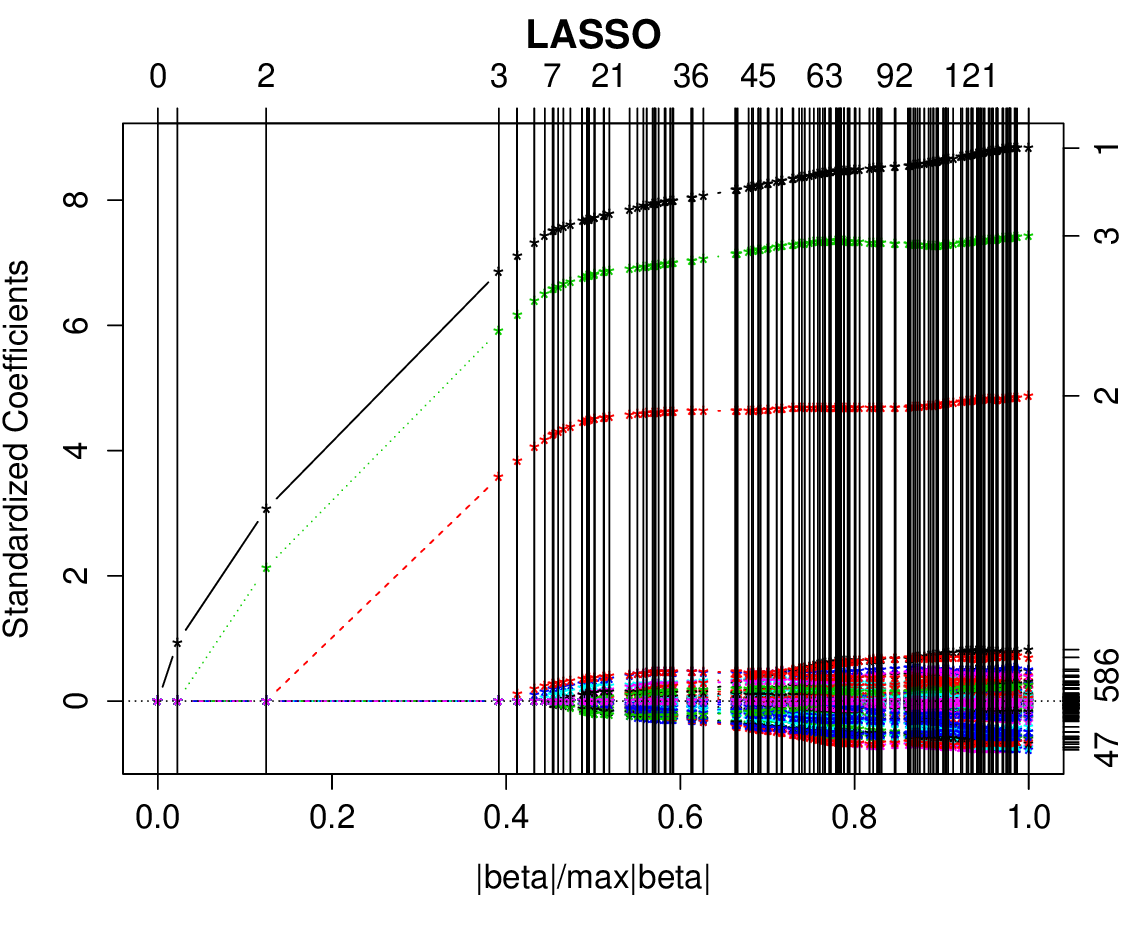}}
	\caption{An example to illustrate the lasso's behavior in the high dimensional setting with different assumptions of noise terms. It reflects that in a situation with standard data and strong sparsity, lasso always chooses the right model no matter the distribution of noise terms.}
	\label{fig:2}
\end{figure}

The results are depicted in Figures \ref{fig:2}. It reflects that in a situation with standard data and strong sparsity, lasso always chooses the right model no matter the distribution of noise terms.

\Section{Proof of Theorem~\ref{thm1}}

Review the lasso estimator
\[
\widehat \beta_n(\lambda_n)\in  \mathop {\text{argmin}}\limits_{\beta\in \mR^{p_n} }\{\dfrac{1}{2}||y_n-X_n\beta||^2_2+\lambda_n ||\beta||\}.
\]

Let $\widehat u_n= \sqrt{n} (\widehat \beta_n -\beta_n)$ and
\[
F_n( \beta_n) =   \dfrac{1}{2}||y_n-X_n\beta_n||^2_2 + \lam_{n} ||\beta_{n}||.
\]

Define $V_n(\widehat u_n)= F_n (\widehat \beta_n)- F_n (\beta_n)$, $C_n = \dfrac{1}{n}X'_n X_n$ and $W_n = \dfrac{1}{\sqrt{n}} X'_n \e_n$. $V_n(\widehat u_n)$ can be written as
\begin{align}
V_n(\widehat u_n)  = \dfrac{1}{2} \widehat u_n' C_n \widehat u_n - \widehat u_n' W_n + \lam_{n} \left( ||\beta_n + \dfrac{\widehat u_n}{\sqrt{n}} || -||\beta_n|| \right). \nonumber
\end{align}

Let $\widehat \beta^{(1)}_n$, $\widehat \beta^{(2)}_n$ and $W^{(1)}_n$, $W^{(2)}_n$ be the first $q$ and last $p_n-q$ elements of $\widehat \beta_n$ and $W_n$ respectively. Similarly, $\widehat u_n^{(1)}$ and $\widehat u_n^{(2)}$ denote the first $q$ and last $p_n-q$ elements of $\widehat u_n$.

Due to $ \big\{ ||W_n||_{\infty}  \leqslant \rho \dfrac{\lam_n}{\sqrt{n}}\big\}$, by (C2) we have
\begin{align}
V_n(\widehat u_n) & \geqslant \dfrac{1-\delta}{2} \left[  \widehat u'^{(1)}_n C_{11,n} \widehat u^{(1)}_n +   \widehat u'^{(2)}_n C_{22,n} \widehat u^{(2)}_n \right]  - \widehat u'_n W_n\notag\\
& -\dfrac{\lam_n}{\sqrt{n}}||\widehat u^{(1)}_n|| + \sum_{j=q+1}^{p_n} |\widehat u_{j,n}|\left (  \dfrac{\lam_n}{\sqrt{n}} - |W_{j,n}|\right ). \nonumber
\end{align}

Since $\widehat u'^{(2)}_n C_{22,n} \widehat u^{(2)}_n \geqslant 0$ and $ \Lambda_{min}(C_{11,n}) $ denotes the smallest eigenvalue of $C_{11,n}$, we have
\begin{align}
V_n(\widehat u_n) & \geqslant ||\widehat u'^{(1)}_n||_2\left\{  \frac{1-\delta}{2}\Lambda_{min}(C_{11,n})  ||\widehat u_n^{(1)}||_2 -||W^{(1)}_n||_2 - \dfrac{\lam_n}{\sqrt{n}}\sqrt{q}\right\} \notag\\
&+ \sum_{j=q+1}^{p_n} |\widehat u_{j,n}|\left (  \dfrac{\lam_n}{\sqrt{n}} - |W_{j,n}|\right )\notag\\
& \geqslant ||\widehat u'^{(1)}_n||_2\left\{ \frac{1-\delta}{2} \Lambda_{min}(C_{11,n}) ||\widehat u_n^{(1)}||_2 -\dfrac{\lam_n}{\sqrt{n}}\sqrt{q} (1+\rho)\right\}  \notag\\
&  + \dfrac{\lam_n}{\sqrt{n}}  (1-\rho)\cdot\sum_{j=q+1}^{p_n} |\widehat u_{j,n}|, \nonumber
\end{align}
Define
\[
M_n \equiv \dfrac{2}{1-\delta} \cdot  \dfrac{\lam_n \sqrt{q}}{\sqrt{n}} (1+\rho)\cdot\Lambda^{-1}_{min}(C_{11,n}).
\]
Then $V_n (\widehat u_n) >0$ depends on
\[
\left \{  || \widehat u_n^{(1)}||_2  >  M_n \right \}.
\]

Since $V_n(0)= 0$, the minimum of $V_n(\widehat u_n)$ cannot be attained at $||\widehat u^{(1)}_n ||_2 > M_n$. Then assume that $\{ \widehat u_n \in \mR^{p_n}:||\widehat u^{(1)}_n ||_2\leqslant M_n, \ \widehat u^{(2)}_n\neq 0 \}$ and (C1) holds. Set $ e_j $ to be the unit vector in the direction of $ j $-th coordinate. Then following inequality holds uniformly:
\begin{align}
\label{eq:V2}
V_n(\widehat u_n) - V_n( \widehat u_n^{(1)}, 0) & =  (\widehat u_n^{(1)})' C_{12,n} \widehat u_n^{(2)} + \frac{1}{2} (\widehat u_n^{(2)})' C_{22,n} \widehat u_n^{(2)}  \notag\\
&+ \dfrac{\lam_{n}}{\sqrt{n}}||\widehat u_n^{(2)}||- (\widehat u_n^{(2)})'W^{(2)}_n \notag\\
& > \sum\limits_{j=q+1}^{p_n}|\widehat u_{j,n}| \left[  \dfrac{\lam_n}{\sqrt{n}} -|W_{j,n}| -\left|\left( (\widehat u_n^{(1)})'C_{12,n}\right)_j\right|  \right] \notag\\
& \geqslant \sum\limits_{j=q+1}^{p_n}|\widehat u_{j,n}| \left[  \dfrac{\lam_n}{\sqrt{n}} (1-\rho) -M_n ||C_{12,n}e_j||_2  \right] \notag\\
& > 0.
\end{align}

Set $ \eta>0 $ such that $ 1-\eta = \dfrac{1-\rho}{1+\rho} \cdot \dfrac{1-\delta}{2} \Lambda_{min}q^{-\frac{1}{2}}$. By (C1) the last inequality of \eqref{eq:V2} holds. Then the minimum of $V_n(u_n)$ can not be attained at $u_n^{(2)} \neq 0$ too, hence we have
\[
\mathop {argmin }\limits_{\widehat u_n \in \mR^{p_n}} V_n(\widehat u_n) \in  \left\{ u_n \in \mR^{p_n} :  ||\widehat u_n^{(1)}||_{2} \leqslant M_n, \widehat u_n^{(2)} = 0  \right\}.
\]

After discussing the model selection consistency of $ \widehat \beta^{(2)}_n $, we now consider about the model selection consistency of $ \widehat \beta^{(1)}_n $. According to the definition of $ \widehat u_n $ and the solution of the lasso, if we want $ \widehat \beta^{(1)}_n \neq 0$ or $ sign(\widehat \beta^{(1)}_n) = sign(\widehat \beta^{(1)}_n) $, the following hold
\[
C_{11,n} \cdot \widehat u^{(1)}_{n} - W^{(1)}_{n} = -\dfrac{\lam_n}{\sqrt{n}} sign(\beta_{n}^{(1)}),
\]
\[
|\widehat u^{(1)}_n/\sqrt{n}| < |\beta^{(1)}_{n}|.
\]

Combining above two restraints of $ \widehat u^{(1)}_n $, the existence of such $ \widehat u^{(1)}_n $ is implied by
\[
|C^{-1}_{11,n}W^{(1)}_{n}| < \sqrt{n}  \left(|\beta^{(1)}_n|  -\dfrac{\lam_n}{n} |C^{-1}_{11,n} \cdot sign(\beta_{n}^{(1)})|  \right).
\]

Since $ C^{-1}_{11,n}W^{(1)}_{n} = C^{-1}_{11,n} \cdot \dfrac{1}{\sqrt{n}} (X^{(1)}_n)' \e $, considering that $\e_i$ are i.i.d random variables with mean $ 0 $ and variance $\sigma^2$ and
\[
\normm{C^{-1}_{11,n} \cdot\frac{1}{\sqrt{n}}(X^{(1)}_n)'}^2_2 =C^{-1}_{11,n},
\]
we have
\[
P([C^{-1}_{11,n}W^{(1)}_{n}]_j > t) \leqslant P([W^{(1)}_{n}]_j > t \cdot  \Lambda_{min}),
\]
where $ \Lambda_{min} = \Lambda_{min}(C_{11,n}) $. Besides, we also have
\[
|C_{11,n}^{-1}sign(\beta_{n}^{(1)})| \leqslant ||C^{-1}_{11,n}||_2 \cdot || sign(\beta^{(1)}_n) ||_2 \leqslant \sqrt{q} \cdot
\Lambda^{-1}_{min}.
\]

By Bonferroni's inequality, we know that if we want to prove
\[
P( \forall j \in S_n, \widehat \beta_{j,n} =0 )	\rightarrow \te^{-nt} \ \text{as} \ n \rightarrow \infty,
\]
it suffices to show that for every $ j \in S_n $,
\[
P(\widehat \beta_{j,n} =0 )	\rightarrow \te^{-nt} \ \text{as} \ n \rightarrow \infty.
\]

Hence we have
\begin{align}
P(\widehat \beta_{j,n} =0 )
& \leqslant P\Big([|W^{(1)}_{n}|]_j \geqslant \sqrt{n} \cdot  \Lambda_{min}
\big(|\beta_{j,n}|  -\dfrac{\lam_n}{n} \Lambda_{min}^{-1}  \sqrt{q}\big)\Big).\nonumber
\end{align}

Let $ G_n = \min \Big\{ \sqrt{n} \cdot  \Lambda_{min} \big(|\beta_{j,n}|  -\dfrac{\lam_n}{n}  \Lambda_{min}^{-1}  \sqrt{q}\big), \
\rho \lam_n / \sqrt{n}  \Big\} $. Then we have
\[
P(\widehat S_n = S_n) \geqslant P(||W_n||_{\infty} \leqslant G_n),
\]
which complets the proof.


\begin{thebibliography}{aa}

\bibitem{bennet1962probability}
Bennet, G.,
\newblock Probability inequalities for sums of independent random variables.
\newblock {\em Journal of the American Statistical Association}, \textbf{57}, 33--45,
1962.

\bibitem{bickel2009simultaneous}
Bickel, P.~J.,  Ritov, Y. and Tsybakov, A.~B.,
\newblock Simultaneous analysis of lasso and dantzig selector.r
\newblock {\em Annals of Statistics}, \textbf{37(4)}, 1705--1732, 2009.

\bibitem{buhlmann2011high}
Buhlmann, P. and Van~de~Geer, S.,
\newblock {\em Statistics for High-Dimensional Data, Methods, Theory and	Applications}.
\newblock Springer, Heidelberg, 2011.

\bibitem{candes2007dantzig}
Candes, E., and Tao, T.,
\newblock The dantzig selector: Statistical estimation when $p$ is much larger
than $n$.
\newblock {\em Annals of Statistics}, \textbf{35(6)}, 2313--2351, 2007.

\bibitem{Efron}
Efron, B., Hastie, T., Johnstone, L. and Tibshirani, R.,
\newblock Least angle regression.
\newblock {\em Annals of Statistics}, \textbf{32(2)}, 407--451, 2004.

\bibitem{fan2014adaptive}
Fan, J.~Q. , Fan, Y.~Y. and Barut, E.,
\newblock Adaptive robust variable selection.
\newblock {\em Annals of Statistics}, \textbf{42(1)}, 324--351, 2014.

\bibitem{fan2001variable}
Fan, J.~Q. and Li, R.~Z.,
\newblock Variable selection via nonconcave penalized likelihood and its oracle
properties.
\newblock {\em Journal of the American Statistical Association},
\textbf{96(456)}, 1348--1360, 2001.

\bibitem{Fan2004}
Fan, J.~Q. and Peng, H.,
\newblock Nonconcave penalized likelihood with a diverging number of
parameters.
\newblock {\em Annals of Statistics}, \textbf{32(3)}, 928--961, 2004.

\bibitem{huang2008asymptotic}
Huang, J., Horowitz, J.~L. and Ma, S.,
\newblock Asymptotic properties of bridge estimators in sparse high-dimensional
regression models.
\newblock {\em Annals of Statistics}, \textbf{36(2)}, 587--613, 2008.

\bibitem{lv2009unified}
Lv, J. C. and Fan, Y.~Y.,
\newblock A unified approach to model selection and sparse recovery using
regularized least squares.
\newblock {\em Annals of Statistics}, \textbf{37(6)}, 3498--3528, 2009.

\bibitem{meinshausen2006high}
Meinshausen, N. and B{\"u}hlmann, P.,
\newblock High-dimensional graphs and variable selection with the lasso.
\newblock {\em Annals of Statistics}, \textbf{34(3)}, 1436--1462, 2006.

\bibitem{meinshausen2010stability}
Meinshausen, N. and B{\"u}hlmann, P.,
\newblock Stability selection.
\newblock {\em Journal of the Royal Statistical Society: Series B (Statistical
	Methodology)}, \textbf{72(4)}, 417--473, 2010.

\bibitem{meinshausen2009lasso}
Meinshausen, N. and Yu, B.,
\newblock Lasso-type recovery of sparse representations for high-dimensional
data.
\newblock {\em Annals of Statistics}, \textbf{37(1)}, 246--270, 2009.

\bibitem{Tibshirani1996}
Tibshirani, R.,
\newblock Regression shrinkage and selection via the lasso.
\newblock {\em Journal of the Royal Statistical Society: Series B.},
\textbf{58}, 267--288, 1996.

\bibitem{tib201201}
Tibshirani, R.~J.,
\newblock The lasso problem and uniqueness.
\newblock {\em Electronic Journal of Statistics}, \textbf{7}, 1456--1490, 2013.

\bibitem{wainwright2009sharp}
Wainwright, M.~J.,
\newblock Sharp thresholds for noisy and high-dimensional recovery of sparsity
using $l_1$-constrained quadratic programming (lasso).
\newblock {\em IEEE Transactions on Information Theory}, \textbf{55(5)}, 2183--2202,
2009.

\bibitem{yyh(elastic)}
Wu, L. and Yang, Y.~H.,
\newblock Nonnegative elastic net and application in index tracking.
\newblock {\em Applied Mathematics and Computation}, \textbf{227}, 541--552, 2014.

\bibitem{yyh(lasso)}
Wu, L., Yang, Y.~H. and Liu, H.~Z.,
\newblock Nonnegative-lasso and application in index tracking.
\newblock {\em Computational Statistics $\&$ Data Analysis}, \textbf{70}, 116--126, 2014.

\bibitem{yyh(adaptive)}
Yang, Y.~H. and Wu, L.,
\newblock Nonnegative adaptive lasso for ultra-high dimensional regression
models and a two-stage method applied in financial modeling.
\newblock {\em Journal of Statistical Planning and Inference}, \textbf{174}, 52--67,
2016.

\bibitem{yuan2006model}
Yuan, M. and Lin, Y.,
\newblock Model selection and estimation in regression with grouped variables.
\newblock {\em Journal of the Royal Statistical Society: Series B}, \textbf{68}, 49--67,
2006.

\bibitem{zhang2010nearly}
Zhang, C.~H.,
\newblock Nearly unbiased variable selection under minimax concave penalty.
\newblock {\em Annals of Statistics}, \textbf{38(2)}, 894--942, 2010.

\bibitem{YuBin06(lasso)}
Zhao, P. and Yu, B.,
\newblock On model selection consistency of lasso.
\newblock {\em Journal of Machine Learning Research}, \textbf{7}, 2541--2563, 2006.

\bibitem{zou2006adaptive}
Zou, H.,
\newblock The adaptive lasso and its oracle properties.
\newblock {\em Journal of the American Statistical Association},
\textbf{101}, 1418--1429, 2006.

\bibitem{zou2005elastic}
Zou, H. and Hastie, T.,
\newblock Regularization and variable selection via the elastic net.
\newblock {\em Journal of the Royal Statistical Society: Series B},
\textbf{67}, 301--320, 2005.

\bibitem{zou2009adaptive}
Zou, H. and Zhang, H.~L.,
\newblock On the adaptive elastic-net with a diverging number of parameters.
\newblock {\em Annals of Statistics}, \textbf{37(4)}, 1733--1751, 2009.

\end{thebibliography}
\end{document}